\documentclass{amsart}
\usepackage{amssymb}
\usepackage{amscd}
\DeclareFontFamily{OMS}{cmsy}{%
\fontdimen16\font=3pt
\fontdimen17\font=3pt}
\makeatletter
\renewcommand{\subsection}{\@startsection{subsection}{2}{\z@}%
{\baselineskip}{0.5\baselineskip}{\bfseries}}
\makeatother
\def\dj{d\kern-.30em\raise1.25ex\vbox{\hrule width .3em height .03em}}
\def\Dj{D\rlap{\kern-.70em\raise0.75ex
\vbox{\hrule width .3em height .03em}}}
\def\dM{\mbox{$\smash{\vphantom{d}^{\raise0.4ex\hbox{$\scriptscriptstyle M$}}%
\!\!\!d}$}}
\def\bla#1{$(${\it #1\/{}}$)$}
\def\U{\mathrm{U}}
\def\SU{\mathrm{SU}}
\def\restr{\restriction}

\def\cal{\mathcal}
\def\Bbb{\mathbb}
\def\e{\epsilon}
\def\k{\kappa}
\def\ad{\mathrm{ad}}
\def\id{\mathrm{id}}
\def\vh{\mathfrak{vh}}
\def\grten{\mathbin{\widehat{\otimes}}}
\def\gr{\mathrm{gr}}

\def\DP{\mathfrak{der}(P)}
\def\DaP{\mathfrak{der}(\aP)}
\def\vDP{\overrightarrow{\mathfrak{der}}(P)}
\def\vDaP{\overrightarrow{\mathfrak{der}}(\aP)}
\def\hor{\mathfrak{hor}}
\def\WM{\Omega_M}
\def\aWP{\Omega(\aP)}
\def\V{\Bbb{V}}
\def\v{\varkappa}
\def\pahh{p_*}
\def\Tcon{\Lambda}
\def\T#1{\mathrm{T}(#1)}
\def\aG{\widetilde{G}}
\def\aA{\widetilde{\cal{A}}}
\def\SV{S(\V)}
\def\zhP{\mathfrak{zh}_P}
\def\azhP{\mathfrak{zh}[{\aP}]}
\def\horP{\hor_P}
\def\ahorP{\mathfrak{h}[\aP]}
\def\aP{\widetilde{P}}
\def\gaP{Q}
\def\aB{\widetilde{\cal{B}}}
\def\gaB{\cal{C}}
\def\Fh{F^\wedge}
\def\Fah{H}
\def\aF{H}
\def\wFh{\widehat{\Fah}}
\def\gFh{K}
\def\aC{\Upsilon}
\def\aCi{\Upsilon_{\inv}}
\def\adaG{\ad}

\def\rig{\wp}
\def\inv{{i\!\hspace{0.8pt}n\!\hspace{0.6pt}v}}
\def\im{\mathrm{im}}
\def\Sum{{\displaystyle\sum}}
\def\tP{\tau}
\def\taP{\tau}
\def\lt#1{[#1]_1}
\def\rt#1{[#1]_2}
\def\sP{\sigma_M}
\def\saP{\sigma_M}
\def\E{\mathrm{E}}
\def\wH{\widetilde{\cal{H}}}
\def\wpi{\widehat{\pi}} 
\renewcommand{\thepage}{\ifnum\value{page}=1 \else\arabic{page}\fi}
\newtheorem{thm}{Theorem}[section]
\newtheorem{pro}[thm]{Proposition}
\newtheorem{lem}[thm]{Lemma}
\theoremstyle{definition}
\newtheorem{defn}{Definition}
\newenvironment{pf}{\proof[\proofname]}{\endproof}
\begin{document}
\title[quantum principal bundles]{AFFINE STRUCTURES\\
ON QUANTUM PRINCIPAL BUNDLES}
\author{micho {\Dj}UR{\Dj}EVICH}
\address{Instituto de Matematicas,
UNAM, Area de la Investigacion Cientifica,
Circuito Exterior, Ciudad Universitaria, M\'exico DF, cp 04510,
MEXICO}

\email{ micho@matem.unam.mx\newline
\indent\indent http://www.matem.unam.mx/{\~\/}micho}

\begin{abstract}
Quantum affine bundles are quantum principal bundles with affine
quantum structure groups.
A general theory of quantum affine bundles is presented. In particular,
a detailed analysis of differential calculi over these bundles is performed,
including the description of a natural differential calculus
over the structure affine quantum group.
A particular attention is given to the study of the
specific properties of quantum affine connections and several
purely quantum phenomena appearing in the context of quantum affine
bundles. Various interesting constructions are presented. In particular, 
the main ideas are illustrated within the example of the quantum Hopf fibration. 
\end{abstract}
\maketitle

\section{Introduction}

The aim of this paper is to study affine quantum principal bundles. These objects
are noncommutative-geometric \cite{c} analogs of classical affine extensions
of principal bundles \cite{KN}. Quantum principal bundles \cite{d1,d2}
generalize classical principal bundles---quantum groups play the role
of the structure groups, and both the base manifold and the bundle are considered as
general quantum objects. In this paper we shall consider very special 
quantum principal bundles, possessing {\it inhomogeneous} 
quantum structure groups. These groups are constructed by applying a 
general quantum affine extension procedure \cite{w-aff} to 
a given (compact matrix) quantum group.

The motivation for this work comes from a high importance
of the classical affine bundles in various areas of standard differential
geometry and theoretical physics. From geometrical point of view,
the formalism of affine bundles gives a unified and powerful language
to express properties and mutual interrelations of entities associated
to the standard geometrical structures
(as for example Riemannian/spin or symplectic manifolds and torsion structures).
As a basic example coming from theoretical physics, let us mention various
geometric formulations of classical general relativity as a gauge theory
associated to Poincare group. Hence, it looks promising to develop a quantum 
version of the theory of affine bundles. In particular, it is interesting to see 
how these objects are related with geometrical structures on quantum spaces, 
described by quantum frame bundles \cite{d-frm2}. The formalism developed here 
may be used in constructing a noncommutative-geometric version of the gauge theory
of gravity. 

The paper is structured in the following way.
In the next section, we shall start from
a quantum principal bundle $P$, with the structure group $G$ and over a given
quantum space $M$. We shall then construct in a natural way another quantum
principal bundle $\aP$ over the same quantum space $M$,
with the structure group $\aG$, which
is the affine extension \cite{w-aff} of $G$ relative
to a given bicovariant bimodule
$\Psi$ over $G$. The left-invariant part of $\Psi$ represents `translational'
degrees of freedom.
After presenting the construction of $\aP$, we shall consider
questions related to differential calculus on affine bundles,
following a general constructive approach to differential calculus
developed in \cite{d-diff}. Our
starting point will be the appropriate graded *-algebra $\horP$, playing the
role of the horizontal forms on the bundle $P$. This algebra will be further naturally
extended to an affine horizontal forms algebra $\ahorP$,
representing horizontal forms on $\aP$. We shall then consider abstract
`covariant derivative' operators
acting in both horizontal forms algebras, and study their mutual relations.
Starting from such operators, it is possible to construct, in a natural manner, a differential
calculus on the affine bundle $\aP$, as well as the associated
differential (bicovariant *-) calculus $\aC$ on $\aG$. A special attention will be given
to the study of the interrelations between differential structures on $P$ and
$\aP$, and to the explicit description of the mentioned calculus
$\aC$ over $\aG$
in terms of the calculus $\Gamma$ over $G$ and the initial bimodule $\Psi$.
As we shall see, an interesting purely quantum phenomena
appears---the left-invariant
part of $\aC$ is not a simple direct sum of the left-invariants
$\Psi_{\inv}$ and $\Gamma_{\inv}$. Instead, it will contain polynomial
components over $\Psi_{\inv}$, generally with arbitrary high degrees.
This is a consequence of the fact that bicovariant bimodules are inherently
braided. We shall also analyze the structure of affine connections,
and the associated curvature, in light of the relations between the calculi
on $P$ and its affine extension $\aP$.

In the last section, concluding remarks and examples are made. At first,
we shall consider quantum frame bundles, and study the specific
properties of their affine extensions. Quantum frame bundles \cite{d-frm2} are
noncommutative-geometric counterparts of the classical frame bundles \cite{KN},
incorporating into the quantum context the concept of a geometrical structure
on a space (like metrics, spinor, symplectic and complex), via the system of 
special `coordinate $1$-forms'. As we shall see, affine extensions of quantum
frame bundles always include (as a subalgebra) a full differential calculus on the
associated quantum plane. 

As a simple but very instructive example, we shall illustrate all the basic elements of the theory 
on the quantum Hopf fibration---which is a
quantum $\U(1)$-bundle over a quantum sphere \cite{p}, the total space $P$ of
which is a quantum $\SU(2)$ group \cite{w-su2}.  If we define $\Psi$ to be a
natural $2$-dimensional bicovariant bimodule over $\U(1)$, coming from the $3D$-calculus
\cite{w-su2} on the quantum $\SU(2)$ group, then the affine extension turns out to be a deformation of the standard
spin 2-covering of $\E(2)$. The corresponding first-order calculus $\aC$ will be 
computed explicitly. Surprisingly, in the truly quantum case this calculus is 
infinite-dimensional. 

We shall also point out some analogies between the formalism of affine
bundles and the formulation of differential calculus \cite{d2} for general
quantum principal bundles.

\section{Affine Quantum Principal Bundles}

This section is devoted to the construction of
quantum principal affine bundles, naturally associated to a given
quantum principal bundle $P$. The structure group
of these bundles is the quantum \cite {w-aff} affine extension $\aG$ of
the structure quantum group $G$ for $P$, with
respect to a given bicovariant bimodule $\Psi$.

\subsection{The Level of Groups}
At first, let us recall the construction of $\aG$. Conceptually, we
follow \cite{w-aff}. Let us consider the tensor bundle algebra
$\T{\Psi}=\Psi^\otimes$.
The coproduct map admits a natural extension to a coassociative
unital homomorphism $\phi\colon\T{\Psi}\rightarrow\T{\Psi}\otimes\T{\Psi}$.
This extension is specified by
$$
\phi(\vartheta)=\ell_\Psi(\vartheta)+\rig_\Psi(\vartheta)
$$
where $\ell_\Psi\colon\Psi\rightarrow\cal{A}\otimes\Psi$ and
$\rig_\Psi\colon\Psi\rightarrow\Psi\otimes\cal{A}$ are the left and the
right action maps respectively. The extended counit
$\e\colon\T{\Psi}\rightarrow\Bbb{C}$ is given by
setting $\e(\vartheta)=0$, for each $\vartheta\in\Psi$. Finally, the
antipode map $\k\colon\cal{A}\rightarrow\cal{A}$
admits a unique antimultiplicative extension
$\k\colon\T{\Psi}\rightarrow\T{\Psi}$, such that
\begin{equation*}
\begin{CD} \Psi @>>> \cal{A}\otimes\Psi\otimes\cal{A}\\
@V{\mbox{$-\k$}}VV @VV{\mbox{$\k\otimes\id\otimes\k$}}V\\
\Psi @<<< \cal{A}\otimes\Psi\otimes\cal{A}
\end{CD}
\end{equation*}
where the horizontal arrows represent the corresponding
twofold coaction and bimodule multiplication respectively.

It is easy to see that $\T{\Psi}$, endowed with the constructed maps,
becomes a Hopf algebra. Furthermore, if $\Psi$ is in
addition a *-covariant bimodule then we have the associated *-involution
$*\colon\Psi\rightarrow\Psi$ such that
$\ell_\Psi$ and $\rig_\Psi$ are hermitian.  The algebra
$\T{\Psi}$ is equipped with the induced *-structure, and
in particular the extended
coproduct map $\phi\colon\T{\Psi}\rightarrow\T{\Psi}\otimes\T{\Psi}$
will be hermitian. The hermicity of the coproduct
implies that $\e$ is hermitian and that the composition map $*\k$
is involutive. In other words, $\T{\Psi}$ is a Hopf *-algebra. In what follows we shall assume that $\Psi$ is equipped
with a *-structure.

Let us denote by $\Psi_{\inv}=\V$ the left-invariant part of
$\Psi$. This space is equipped with a natural right $\cal{A}$-comodule
structure $\v\colon\V\rightarrow\V\otimes\cal{A}$ and with a right
$\cal{A}$-module structure $\circ\colon\V\otimes\cal{A}\rightarrow\V$.
Here $\v$ is the restriction of the right-action on $\V$, and $\circ$ is
given by the formula $$\{\,\}\circ a=\k(a^{(1)})\{\,\}a^{(2)}.$$
Furthermore, the space $\V$ is *-invariant, and the following compatibility
relations hold
\begin{gather*}
\v(\theta\circ a)=\sum_k(\theta_k\circ a^{(2)})\otimes
\k(a^{(1)})c_ka^{(3)}\\
\v*=(*\otimes *)\v\qquad (\theta\circ a)^*=\theta^*\circ\k(a)^*.
\end{gather*}
where $\v(\theta)=\Sum_k\theta_k\otimes c_k$. All these maps
admit natural extensions to the tensor algebra $\V^\otimes$,
and the same compatibility relations hold for the extended maps (which
for simplicity will be denoted by the same symbols). The extensions are
fixed by
\begin{gather*}
\v(\theta\eta)=\v(\theta)\v(\eta)\qquad (\theta\eta)^*=\eta^*\theta^*\\
(\theta\eta)\circ a=(\theta\circ a^{(1)})(\eta\circ a^{(2)})\qquad
1\circ a=\e(a) 1.
\end{gather*}

The tensor algebra $\V^\otimes\subseteq\T{\Psi}$ is right $\phi$-invariant, in the sense that
$$
\phi[\V^\otimes]\subseteq\V^\otimes\otimes\T{\Psi}.
$$

Let $\tau\colon\Psi\otimes_{\cal{A}}\Psi\rightarrow\Psi\otimes_{\cal{A}}
\Psi$ be the braid operator, naturally associated \cite{w2} to the bicovariant
bimodule $\Psi$. This operator is
left/right covariant, and in particular it is completely
determined by its restriction
$$
\tau\colon\V^{\otimes 2}\rightarrow\V^{\otimes 2},\qquad
\tau(\eta\otimes\theta)=\sum_k\theta_k\otimes(\eta\circ c_k),
$$
where $\Sum_k\theta_k\otimes c_k=\v(\theta)$.

Now let us consider the $\tau$-symmetrizer maps
$Y_k\colon\T{\Psi}^k\rightarrow\T{\Psi}^k$, defined by
$$
Y_k=\sum_{\pi\in S_k} \tau_\pi,
$$
where $\tau_\pi$ is the operator obtained by replacing the
transpositions figuring in any minimal decomposition of the permutation
$\pi$ by the corresponding $\tau$-twists (this definition is consistent,
due to the braid equation---in other words $\tau_\pi$ does not depend of
the choice of minimal decompositions).
The following factorizations hold:
$$
Y_{k+l}=Y_{kl}(Y_k\otimes Y_l)=(Y_k\otimes Y_l)M_{kl},
$$
where
$$ Y_{kl}=\sum_{\pi\in S_{kl}} \tau_\pi \qquad
M_{kl}=\sum_{\pi\in S_{kl}} \tau_{\pi^{-1}}, $$
and $S_{kl}\subset S_{k+l}$ consists of permutations preserving the
orders of the first $k$ and the last $l$ factors.

In particular, it follows that the space $\ker(Y)$ is a two-sided ideal in
$\T{\Psi}$, where $Y\colon\T{\Psi}\rightarrow\T{\Psi}$ is the corresponding
`total symmetrizer'. Furthermore, the following commutation relations hold:
$$
*\tau*=\tau^{-1}\qquad \k\tau=\tau\k.
$$
This implies that $\ker(Y)^*=\ker(Y)$ and
$\k\bigl[\ker(Y)\bigr]=\ker(Y)$.

We have a natural identification $\Psi\leftrightarrow\cal{A}\otimes\V$,
and in particular $\T{\Psi}\leftrightarrow\cal{A}\otimes\V^\otimes$.
In terms of this identification, the braid $\tau$ and the
extended coproduct $\phi$ are connected by the formula
\begin{equation}\label{ext-coprod}
\phi(a\otimes\vartheta)=\phi(a)\sum_{k+l=n} (\v_k \otimes\id^l)
M_{kl}(\vartheta).
\end{equation}
It follows that
$$ \phi\bigl[\ker(Y)\bigr]\subseteq\ker(Y)\otimes\T{\Psi}+
\T{\Psi}\otimes\ker(Y).$$
In other words, $\ker(Y)\subset\T{\Psi}$ is a Hopf *-ideal and
therefore we can pass to the factor Hopf *-algebra
$$\aA=\T{\Psi}/\ker(Y).$$
We shall denote by the same symbols $\phi,\e,\k,*$ the corresponding
projected maps.

The ideal $\ker(Y)$ is bicovariant. In particular, $\aA$ is
bicovariant and so we have the following natural decompositions:
$$ \ker(Y)=\cal{A}\otimes\ker\bigl(Y{\restriction}\Psi_{\inv}^\otimes
\bigr)\qquad \aA\leftrightarrow\cal{A}\otimes \SV,$$
with an independent definition
\begin{gather*}
\SV=\Psi_{\inv}^\otimes/
\ker\bigl(Y{\restriction}\Psi_{\inv}^\otimes\bigr)\\
\SV\leftrightarrow\im(Y)\hookrightarrow\V^\otimes.
\end{gather*}

We shall also denote by the symbols $\v$ and $\circ$ the corresponding
right coaction and the right $\cal{A}$-module structure on
$\SV$ respectively. We have
$$ \phi[\SV]\subseteq\SV\otimes\aA,$$
and we shall denote by the symbol
$\widehat{\v}\colon\SV\rightarrow\SV\otimes\aA$ the corresponding
restriction map.

\subsection{The Level of Bundles}

Let $M$ be a quantum space and $P=(\cal{B},i,F)$ a quantum principal 
$G$-bundle \cite{d2} over $M$. Here $\cal{B}$ is a *-algebra representing the
quantum space $P$, and $F\colon\cal{B}\rightarrow\cal{B}\otimes\cal{A}$
is a counital coassociative *-homomorphism playing the role of the
dualized right action of $G$ on $P$. Furthermore
$i\colon\cal{V}\rightarrow\cal{B}$ is a *-isomorphism between
a *-algebra $\cal{V}$ representing the quantum space $M$ and the $F$-fixed
point subalgebra of $\cal{B}$. Therefore we can identify $\cal{V}$ with
its image in $\cal{B}$. Finally, we require that the action $F$ is
free, in the sense that a map $X\colon\cal{B}\otimes_{\cal{V}}\cal{B}
\rightarrow\cal{B}\otimes\cal{A}$ given by $X(q\otimes b)=qF(b)$ is
surjective.

We are now going to extend the bundle $P$, by including
`translational' degrees of freedom. This will give us a quantum principal
$\aG$-bundle $\aP$, over the same quantum space $M$.

At first, let us define a *-algebra $\aB$ as
\begin{equation}
\aB=\cal{B}\otimes\SV,
\end{equation}
at the level of vector spaces, and let us assume that $\aB$ is
equipped with a *-algebra
structure specified by cross-product type formulae
\begin{align}
(q\otimes\theta)(b\otimes\eta)&=
\sum_kqb_k\otimes(\theta{\circ} c_k)\eta\\
(b\otimes\theta)^*&=\sum_k b_k^*\otimes(\theta^*{\circ} c_k^*),
\end{align}
where $F(b)=\Sum_kb_k\otimes c_k$. We see that
both $\cal{B}$ and $\SV$ are understandable as *-subalgebras of $\aB$.
The algebra $\aB$ is naturally graded, with the grading induced from $\SV$.

The maps $F$ and $\widehat{\v}$ admit a natural common extension to a
*-homomorphism $\aF\colon\aB\rightarrow\aB\otimes\aA$.
More precisely $\aF$ is the direct product of
actions $F$ and $\widehat{\v}$, and in particular
\begin{equation}
(\aF\otimes\id)\aF=(\id\otimes\phi)\aF\qquad
(\id\otimes\e)\aF=\id.
\end{equation}

To prove that $\aF$ is really a *-homomorphism, it is sufficient to
check that the commutation relations between $\cal{B}$ and $\SV$ are
preserved. We have
\begin{equation*}
\begin{split}
\vartheta b&=\sum_k b_k(\vartheta\circ c_k)\longrightarrow
\sum_{kl}(b_k\otimes c_k^{(1)})\Bigl\{(\vartheta_l\circ c_k^{(3)})\otimes
\k(c_k^{(2)})d_l\k(c_k^{(4)})\Bigr\}\\
&=\sum_{kl}b_k(\vartheta_l\circ c_k^{(1)})\otimes d_lc_k^{(2)}=
\sum_{kl}(\vartheta_l\otimes d_l)(b_k\otimes c_k)
=\aF(\vartheta)\aF(b),
\end{split}
\end{equation*}
where $\varphi\in\cal{B}$ and $\vartheta\in\SV$ with
$\Sum_l\vartheta_l\otimes d_l=\widehat{\v}(\vartheta)$. In other words,
$\aG$ acts by `automorphisms' on $\aB$.

\begin{lem}
The fixed-point subalgebra for the total action $\aF$ on
$\aB$ coincides with the algebra $\cal{V}\leftrightarrow\cal{V}\otimes\Bbb{C}$.
\end{lem}

\begin{pf}
Let us consider an arbitrary element
$q=\Sum_{b\vartheta}b\otimes\vartheta$, with linearly
independent and homogeneous $\vartheta\in\SV$ and arbitrary $b\in\cal{B}$.
From the form of the action $\widehat{\v}$, it follows that
$$
\aF (q)=\sideset{}{^*}\sum_{b\vartheta}F(b)\vartheta+
q'
$$
where the summation is performed over pairs $(b,\vartheta)$ having
the maximal degree, and $q'\in\aB\otimes\aA$ consists of
summands with lower degrees (here we have compared degrees relative
to the second tensoriand). It follows that
\begin{equation}
\cal{B}=\aF^{-1}\Bigl\{\aB\otimes\cal{A}\Bigr\}.
\end{equation}
In particular $q$ is $\aF$-invariant iff
$q=b\otimes 1\leftrightarrow b$, and $b$ is $F$-invariant.
\end{pf}

As already mentioned, the freeness of the action $F$ is expressed as
the surjectivity of the natural map
$X\colon\cal{B}\otimes_{\cal{V}}\cal{B}\rightarrow\cal{B}\otimes\cal{A}$.
As explained in \cite{d-ext}, the map $X$ is also injective for
compact structure groups $G$, 
so that $P$ defines a Hopf-Galois extension \cite{s-ext}. 
In particular, we can introduce the quantum `translation
map' by inverting $X$ over $\cal{A}$---in other words
we define $\tP\colon\cal{A}\rightarrow\cal{B}\otimes_{\cal{V}}\cal{B}$ by
$$
\tP(a)=X^{-1}(1\otimes a),
$$
so that
$$
b\otimes a=X\bigl[b\tP(a)\bigr] \quad\qquad \forall\, b\in\cal{B}, \,a\in\cal{A}.
$$

If $G$ is a general non-compact structure group then the existence of
$\tP$ is a non-trivial condition on the bundle $P$. However, the main constructions
of this paper will still work if we assume that $\tP$ exists (that is, if $P$ is a
Hopf-Galois extension---equivalent to the bijectivity of $X$).

Now we shall prove that $\tP$ admits a natural extension
$\taP\colon\aA\rightarrow\aB\otimes_{\cal{V}}\aB$, which is
actually the affine translation map. In what follows, we shall also use the 
symbolic notation $\tP(\,)=\lt{\,}\otimes\rt{\,}$.

We shall start by writing the explicit formula for the extended translation
map. Let $\taP\colon\aA\rightarrow\aB\otimes_{\cal{V}}\aB$ be a
linear map defined by
\begin{equation}\label{taP1}
\taP(a\otimes\vartheta)=\sum_{i\geq 0}\sum_{\alpha,\beta\in S[i]}(-1)^i 
\mathrm{I}\{t_{\beta i}\} \lt{c_{\beta\alpha}^i}
\tP(a)\rt{c_{\beta\alpha}^i}\vartheta_{\alpha i}.
\end{equation}
In the above formula $\vartheta\in\SV^n$ and the elements $\bigl\{t_{\alpha i}\mid\alpha\in S[i]\bigr\}$ form a basis in the space $\SV^i$. Furthermore, the map $\mathrm{I}\colon\SV\rightarrow
\SV$ is the canonical `total braided inverse'. More precisely, this map is originally
defined on $\V^\otimes$ as the unital grade-preserving map totally reversing 
the order with the help of the braiding $\tau$---and then it is projected down 
to $\SV$. Finally, we have put \begin{gather*}
\v(t_{\alpha i})=\sum_{\beta\in S[i]} t_{\beta i}\otimes c_{\beta\alpha}^i\\
\sum_{\alpha\in S[i]} t_{\alpha i}\otimes\vartheta_{\alpha i}=M_{i n-i}(\vartheta).
\end{gather*}
The consistency of the above formula follows from the fact that
operators $M_{i n-i}$ are projectable down to $\SV^n$.

In particular, we see that
\begin{equation}\label{taP2}
\taP{\restr}\cal{A}=\underbrace{\tP\colon\cal{A}
\rightarrow\cal{B}\otimes_{\cal{V}}\cal{B}}_{\text{\normalsize\it the original map}}
\qquad\qquad \taP(\theta)=-\sum_k\theta_k\tP(c_k)+ 1\otimes\theta,
\end{equation}
where $\theta\in\V$ and $\Sum_k\theta_k\otimes c_k=\v(\theta)$.

\begin{lem} The extended map satisfies  
\begin{align*}
\lt{\psi}\otimes \aF\rt{\psi}&=\lt{\psi^{(1)}}\otimes\rt{\psi^{(1)}}\otimes \psi^{(2)}\\
\lt{\psi}\rt{\psi}&=\e(\psi)1\\
{}^{\mathrm{op}}\!\aF\lt{\psi}\otimes \rt{\psi}&=\psi^{(1)}\otimes \lt{\psi^{(2)}}\otimes\rt{\psi^{(2)}}
\end{align*}
where ${}^{\mathrm{op}}\!\aF\colon\aB\rightarrow\aA\otimes\aB$ is the opposite
action, defined by ${}^{\mathrm{op}}\!\aF=(\k^{-1}\otimes\id)(\circlearrowright)\aF$. 
\end{lem}

\begin{pf}
These identities follow by direct calculations, applying formula \eqref{ext-coprod}
for the extended coproduct and the antipode formula
\begin{equation}\label{ext-ant}
\k(a\otimes\vartheta)=\Bigl\{
(-1)^{\partial\vartheta}\sum_k\mathrm{I}\{\vartheta_k\}\k(c_k)\Bigr\}\k(a)
\end{equation}
which holds both in $\aA$ and $\T{\Psi}$. Finally, let us observe that
the original translation map $\tP\colon\cal{A}\rightarrow\cal{B}
\otimes_{\cal{V}}\cal{B}$ satisfies our identities by construction. 
\end{pf}

The above lemma implies that the extended
$X\colon\aB\otimes_{\cal{V}}\aB\rightarrow\aB\otimes\aA$ is also 
bijective, and that we have $$X^{-1}(1\otimes\psi)=\tP(\psi).$$
In other words $\aB$ equipped with $\aF$ determines a Hopf-Galois 
extension and the extended
$\taP$ is the affine translation map. In particular, the
action $\aF$ is free, and we have a quantum principal $\aG$-bundle
$\aP=(\aB,i,\aF)$ over the quantum space $M$.

\begin{defn} The bundle $\aP$ is called {\it the quantum affine extension}
of the quantum principal bundle $P$.
\end{defn}

We can interpret the quantum space $\aP$ as a bundle over $P$. The
fibering of $\aP$ over $P$ is trivial, and the fiber is the quantum affine space associated
to $\V$ and $\tau$---described by the symmetric algebra $\SV$. The algebra
$\cal{B}$ is the fixed-point subalgebra for the action $\id\otimes\widehat{\v}
\colon\aB\rightarrow\aB\otimes\aA$. Let
$\pahh\colon\aB\rightarrow\cal{B}=\aB^0$ be the canonical projection map.
If we interpret $F$ as the restricted action of $\aG$ we can write
$$
F\pahh=(\pahh\otimes\id)\aF,
$$
and therefore $\aP$ is geometrically interpretable as the extension of
the bundle $P$.

As explained in \cite{d-qpbb}, if a quantum principal bundle
admits the translation map then there exists an intrinsic braid operator,
twisting the functions of the bundle. In our context this means that
we can introduce a braiding $\sP\colon\cal{B}\otimes_{\cal{V}}\cal{B}
\rightarrow\cal{B}\otimes_{\cal{V}}\cal{B}$, and its affine extension
$\saP\colon\aB\otimes_{\cal{V}}\aB\rightarrow\aB\otimes_{\cal{V}}\aB$.

\begin{lem}
There exists an intrinsic braiding
$\saP\colon\aB\otimes_{\cal{V}}\aB\rightarrow\aB\otimes_{\cal{V}}\aB$,
extending the initial braiding $\sP$. We have
\begin{equation}\label{a-braid1}
\saP(\theta\otimes\varphi)=\sum_k[\theta_k,\varphi]\tP(c_k)+\varphi\otimes\theta.
\end{equation}
for each $\varphi\in\aB$ and $\theta\in\V$.
\end{lem}

\begin{pf} According to \cite{d-qpbb}, the braidings are defined by the formula
\begin{equation}\label{bbraid}
\sP(b\otimes \psi)=\sum_k b_k\psi\lt{c_k}\otimes\rt{c_k},
\end{equation}
where $\Sum_k b_k\otimes c_k=\aF(b)$. Taking this into account, formula 
\eqref{a-braid1} directly follows from \eqref{taP2}.
\end{pf}

Let us observe that formula \eqref{a-braid1} uniquely determines
the extended braiding $\saP$.
This follows from the fact that the braiding is always functorial
relative to the product in the bundle algebra \cite{d-qpbb}. Explicitly, the formula 
\eqref{a-braid1} is generalized to the following expression describing the braiding 
between $\SV$ and $\aB$ (the notation is
the same as in the definition of the affine translation map):

\begin{pro} We have
\begin{equation}\label{a-braid2}
\saP(\vartheta\otimes\varphi)=\varphi\otimes\vartheta +\sum_{i>0}\sum_{\alpha,\beta\in S[i]}\{t_{\beta i}\vert\varphi\}_\tau
\lt{c_{\beta\alpha}^i}\otimes\rt{c_{\beta\alpha}^i}\vartheta_{\alpha i},
\end{equation} 
where we have introduced `braided multicommutators' 
$\{\mid\}$ defined by
\begin{equation}
\{\vartheta\vert\varphi\}_\tau=\sum_{i\geq 0}\sum_{\alpha \in S[i]}
(-1)^{n-i}t_{\alpha i}\varphi\mathrm{I}\{\vartheta_{\alpha i}\}. 
\end{equation}
\end{pro}

\begin{pf} 
This identity follows directly from \eqref{bbraid} and \eqref{taP1}.  
It is also possible to derive \eqref{a-braid2} from \eqref{a-braid1}, 
inductively applying the `left variant' of the mentioned functoriality property. 
\end{pf}

\section{Affine Horizontal Forms $\&$ Related Calculus}\label{sec:ahor}

\subsection{Horizontal Forms}

In this subsection we shall refine our construction of the affine bundles,
in order to include the analogs of
horizontal forms in the game. This is the first step
towards constructing the appropriate complete differential calculus on 
affine quantum principal bundles.

Let us consider a quantum principal $G$-bundle
$P=(\cal{B},i,F)$ over a quantum space $M$. Let $\horP$ be
graded *-algebra such that $\cal{B}=\horP^0$, equipped with a coassociative
grade-preserving
*-homomorphism $\Fh\colon\horP\rightarrow\horP\otimes\cal{A}$ extending
the action $F$. Let $\WM\subseteq\horP$ be the corresponding fixed-point
*-subalgebra. Let us assume that $\WM$ is equipped with a
hermitian first-order differential $\dM\colon\WM\rightarrow\WM$. We shall
interpret the elements of $\horP$ as `abstract horizontal forms'.
In the framework of this interpretation, the elements
of $\WM$ naturally correspond to the forms on the base.

Let $\zhP\subseteq\horP$ be the graded commutant of $\WM$ in $\horP$.
This is a graded *-subalgebra of $\horP$, invariant under the action $\Fh$.
As explained in \cite{d-diff}, the formula
$$
\xi\circ a=\lt{a}\xi\rt{a},
$$
consistently defines a right $\cal{A}$-module structure
in the space $\zhP$. This map, together with $\Fh$ and $*$, determines
a bicovariant *-bimodule structure over $G$. The following commutation
relation holds:
\begin{equation}
\xi\varphi=(-1)^{\partial\varphi\partial\xi}\sum_k\varphi_k(\xi\circ c_k),
\end{equation}
where $\varphi\in\horP$.
The algebra $\zhP$ is braided-commutative, in a natural manner.

Let us now prove a useful technical lemma.

\begin{lem} Let $\Delta\colon\horP\rightarrow\horP$ be an arbitrary
graded derivation, intertwining the map $\Fh$. Then $\Delta$ is reduced in
the space $\zhP$, and
\begin{equation}
\Delta(\xi\circ a)=\Delta(\xi)\circ a,
\end{equation}
for each $\xi\in\zhP$ and $a\in\cal{A}$.
\end{lem}

\begin{pf} The fact that $\Delta$ is reduced in $\WM$, together with the
graded Leibniz rule and the definition of $\zhP$, implies that
$\Delta$ is reduced in $\zhP$. We have
\begin{multline*}
\Delta(\xi\circ a)=(\Delta\lt{a})\xi\rt{a}
+(-1)^{\partial\xi\partial\Delta}\lt{a}\xi\Delta\rt{a} +\lt{a}\Delta(\xi)\rt{a}\\=
\Delta\lt{a^{(1)}}\rt{a^{(1)}}(\xi\circ a^{(2)})\\{}+
\lt{a^{(1)}}\Delta\rt{a^{(1)}}(\xi\circ a^{(2)})+\Delta(\xi)\circ a\\
=\Delta\bigl\{\lt{a^{(1)}}\rt{a^{(1)}}\bigr\}(\xi\circ a^{(2)})+\Delta(\xi)\circ a
=\Delta(\xi)\circ a,\\
\end{multline*}
because of the identity $\lt{a}\rt{a}=\e(a)1$.
\end{pf}

Let us denote by $\DP$ the real affine space
of all hermitian $\Fh$-covariant
first-order antiderivations $D\colon\horP\rightarrow\horP$ extending
the differential $\dM\colon\WM\rightarrow\WM$. In what follows, it will be assumed that
$\DP\neq\emptyset$. In accordance with \cite{d-diff}, we adopt the following

\begin{defn}
The elements of $\DP$ are called {\it bundle derivatives} for $P$.
The associated vector  space $\vDP$
consists of hermitian first-order right-covariant antiderivations
$E$ on $\hor_P$ which vanish on $\WM$.
\end{defn}

We shall interpret the elements
of $\DP$ as covariant derivative maps for the bundle $P$. This interpretation
is justified after constructing the appropriate complete calculus on the
bundle, such that bundle derivatives are in $1$--$1$ correspondence with
covariant derivatives of regular connections \cite{d2}. The construction of the
complete calculus on the bundle is presented in full details in \cite{d-diff}. 

As explained in \cite{d-diff},
to each $D\in\DP$ we can intrinsically associate a linear map
$\varrho_D\colon \cal{A}\rightarrow \hor_P$, with the help of the
formula
\begin{equation}\label{d2-curv}
D^2(\varphi)=-\sum_k\varphi_k\varrho_D(c_k),
\end{equation}
where $\Sum_k\varphi_k\otimes c_k=\Fh(\varphi)$. This map plays the
role of {\it the curvature} of $D$. It is uniquely defined by $D$.
Moreover, the following equalities hold
\begin{gather}
\Fh\varrho_D(a)=(\varrho_D\otimes\id)\ad(a)\\
D\varrho_D(a)=0\\
\varrho_D[\kappa(a)^*]=-\varrho_D(a)^*\\
\varrho_D(a)\varphi=\sum_k\varphi_k
\varrho_D(ac_k)\label{varrho-hor}
\end{gather}
for each $a\in\ker(\e)$ and $\varphi\in\horP$. 
 
Furthermore, for each $E\in\vDP$
there exists a unique $\chi_E\colon \cal{A}\rightarrow \horP$
such that
\begin{equation}
E(\varphi)=-(-1)^{\partial\varphi}\sum_k \varphi_k
\chi_E(c_k), \end{equation}
for each $\varphi\in \horP$. The following equalities hold:
\begin{gather}
\Fh\chi_E(a)=(\chi_E\otimes \id)\ad(a)\\
\chi_E[\kappa(a)^*]=-\chi_E(a)^*\\
\chi_E(a)\varphi=(-1)^{\partial\varphi}\sum_k\varphi_k\chi_E
(ac_k)\label{chi-hor}
\end{gather}
for each $a\in \ker(\e)$ and $\varphi\in \horP$.
Let us also observe that 
$$ \varrho_D(1)=\chi_E(1)=0. $$

Applying the construction from the previous section
to the pair $(\horP,\Fh)$ we obtain a *-algebra 
$$\ahorP\leftrightarrow\horP\otimes\SV,$$
where the product and the conjugation are given by the same
cross-product type rules, 
together with the action $\Fah\colon\ahorP\rightarrow\ahorP\otimes\aA$.
We have the following interesting characterization
$$
\horP=\horP\otimes\Bbb{C}=\Fah^{-1}\bigl(\ahorP\otimes\cal{A}\bigr).
$$
The algebra $\ahorP$ is bigraded, in a natural manner.

The elements of the constructed algebra $\ahorP$ play the role of
horizontal differential forms on the quantum affine bundle $\aP$.
Restricted to $\horP$, the action $\Fah$ coincides with $\Fh$. In particular,
the $\Fah$-fixed-point subalgebra of $\ahorP$ coincides with $\WM$.
We shall denote by the same symbol
$\pahh\colon\ahorP\rightarrow\horP$ the canonical projection map.

In a similar way, let us introduce the real
affine space $\DaP$ and its vectorization $\vDaP$.
The maps from $\DaP$ play the role of affine covariant derivatives.
Actually, this interpretation of $\DP$ and $\DaP$
is fully justified after constructing the intrinsic complete
differential calculi over the bundles $P$ and $\aP$. The maps
from $\DP$ and $\DaP$ become `true' covariant derivatives
associated to corresponding regular connections, in the sense of \cite{d2}.

Let us consider an arbitrary $D\in\DaP$. The intertwining property implies
\begin{align*}
D(\horP)&\subseteq\horP\\
D(\V)&\subseteq\horP.
\end{align*}
Moreover, because of the graded Leibniz rule, every such $D$ is completely
determined by its restrictions
$$ (D{\restr}\horP)\in\DP, \qquad\quad
(D{\restr}\V)=\lambda\colon\V\rightarrow\horP.$$

\begin{pro} \bla{i} Conversely, let us consider an arbitrary first-order linear map
$\lambda\colon\V\rightarrow\horP$ satisfying
\begin{align}
\Fh\lambda&=(\lambda\otimes\id)\v\label{v-L}\\
\lambda(\theta)\varphi&=(-1)^{\partial\varphi}\sum_k
\varphi_k\lambda(\theta\circ c_k),\label{hor-L}
\end{align}
where $\varphi\in\horP$ and $\Sum_k\varphi_k\otimes c_k=\Fh(\varphi)$. Then every
$\Fh$-covariant first-order antiderivation
$D\colon\horP\rightarrow\horP$ admits a unique extension
to a first-order
antiderivation $D\colon\ahorP\rightarrow\ahorP$ such that
$D(\theta)=\lambda(\theta)$, for each $\theta\in\V$. 

\smallskip
\bla{ii} This map is
automatically $\Fah$-covariant. It is hermitian, if and only if both $\lambda$ and
the initial $D$ are hermitian.
\end{pro}

\begin{pf} At first, we have to check the compatibility of the definition
of the extended covariant derivative,
with the relations defining the algebra $\ahorP$. A direct
computation gives
\begin{equation*}
\begin{split}
0=\theta\varphi-\sum_k\varphi_k(\theta\circ c_k)\longrightarrow&
\lambda(\theta)\varphi+\theta D(\varphi)\\
{}-&\sum_kD(\varphi_k)(\theta\circ c_k)-(-1)^{\partial\varphi}
\varphi_k\lambda(\theta\circ c_k)=0,
\end{split}
\end{equation*}
because of the $\Fh$-covariance of $D\colon\horP\rightarrow
\horP$ and the definition of $\lambda$. Furthermore,
$$
\SV^n\ni\vartheta\longrightarrow
(\{\}^\wedge\otimes\id^{n-1})M_{1n-1}(\vartheta)
$$
which is consistent because of the invariance of $\ker(Y_n)$ under all
operators $M_{kl}$, where $k+l=n$. Therefore, our definition together
with the graded Leibniz rule consistently and uniquely
determines the extension $D\colon\ahorP\rightarrow\ahorP$. Finally, to
prove \bla{ii} it is sufficient to check that $\V$ is covariant
under the action of $*$ and $\Fah$. The hermicity is trivial, while
$$
(D\otimes\id)\Fah (\theta)=(D\otimes\id)\bigl[1\otimes\theta
+\v(\theta)\bigr]=(\lambda\otimes\id)\v(\theta)=\Fah
D(\theta),
$$
for each $\theta\in \V$. This completes the proof.
\end{pf}

Let us denote by $\Tcon$ the graded space of all maps
$\lambda\colon\V\rightarrow\horP$ intertwining $\v$ and $\Fh$ and satisfying
\begin{equation}\label{hor-L2}
\lambda(\theta)\varphi=(-1)^{\partial\varphi\partial\lambda}\sum_k 
\varphi_k\lambda(\theta{\circ}c_k). 
\end{equation}
In particular, the above equality implies that every map $\lambda\in\Tcon$ takes the 
values from $\zhP$. Furthermore, \eqref{hor-L2} is equivalent to the intertwining property
\begin{equation}
\lambda(\theta\circ a)=\lambda(\theta)\circ a. 
\end{equation}

The formula
$$ \lambda^*(\theta)=\lambda(\theta^*)^*$$
determines a natural conjugation on this space. According to the above analysis,
\begin{equation}\label{dec-aD}
\DaP_{\Bbb{C}}=\DP_{\Bbb{C}}\oplus \Tcon^1,
\end{equation}
which corresponds to the classical decomposition \cite{KN} of the affine connections.

Similarly, we can introduce the affine transition maps
$E\colon\ahorP\rightarrow\ahorP$. They form a real affine space $\vDaP$.
The following natural decomposition holds
\begin{equation}\label{dec-aE}
\vDaP_{\Bbb{C}}=\vDP_{\Bbb{C}}\oplus \Tcon^1,
\end{equation}
in accordance with \eqref{dec-aD} and the definition of $\vDaP$.

We are going to analyze relations between the affine curvature map (defined in the 
framework of $\aP$) and the standard curvature map (defined at the level of $P$).

Let $\varrho_D\colon\aA\rightarrow\ahorP$ be the curvature
map associated to the affine covariant derivative
$D$. It is the extension of the curvature
$\varrho_D\colon\cal{A}\rightarrow\horP$, associated to $D\colon\horP
\rightarrow\horP$. 
 
\begin{pro}\label{pro-affcurv} We have
\begin{equation}\label{ext-r}
\varrho_D(\theta)=-\sum_k\theta_k\varrho_D(c_k)-D\lambda(\theta),
\end{equation}
where $\theta\in \V$ and $\Sum_k\theta_k\otimes c_k=\v(\theta)$.
\end{pro}

\begin{pf} Applying the definitions of $D$, $\varrho_D$ and $\lambda$
we obtain
\begin{equation*}
D^2(\theta)=D\lambda(\theta)=-\varrho_D(\theta)-\sum_k
\theta_k\varrho_D(c_k),
\end{equation*}
in other words \eqref{ext-r} holds.
\end{pf}

Similarly, the maps $E\colon\ahorP\rightarrow\ahorP$ are naturally determined by
linear functionals $\chi_E\colon\aA\rightarrow\ahorP$. These maps extend
previously considered $\chi_E\colon\cal{A}\rightarrow\horP$, and we have
\begin{equation}\label{ext-k}
\chi_E(\theta)=-\sum_k\theta_k\chi_E(c_k)-E(\theta),
\end{equation}
for each $\theta\in\V$.

It is worth noticing that principal relations \eqref{d2-curv}--\eqref{chi-hor} 
are preserved at the affine level.
In particular, formulae \eqref{ext-r}--\eqref{ext-k} completely determine
the affine extensions of the original maps $\varrho_D$ and $\chi_E$.

\subsection{The Induced Affine Calculus}

Let us consider the algebra $\horP$. As explained in \cite{d-diff},
the system of maps
$\varrho_D,\chi_E\colon\cal{A}\rightarrow\horP$ naturally determines a
bicovariant *-calculus $\Gamma$ over $G$. By definition, this calculus
is based on the right $\cal{A}$-ideal $\cal{R}\subset\ker(\e)$
consisting of the elements annihilated by all $\varrho_D$ and $\chi_E$. In order to get
a nontrivial calculus, we shall assume that at least one bundle derivative has a non-zero
curvature. Geometrically, $\Gamma$ is the minimal calculus on $G$ compatible with
the internal structure of horizontal forms $\horP$.

By construction, the maps $\varrho_D$ and $\chi_E$ are factorizable through $\cal{R}$, 
and hence interpretable as $\varrho_D,\chi_E\colon\Gamma_{\inv}\rightarrow \horP$. Let us
also observe that \eqref{varrho-hor}/\eqref{chi-hor} imply that $\varrho_D,\chi_E$ take 
the values from $\zhP$, and that 
\begin{equation}
\varrho_D(\xi\circ a)=\varrho_D(\xi)\circ a\qquad\quad\chi_E(\xi\circ a)=\chi_E(\xi)\circ a. 
\end{equation} 

On the other hand, applied at the level of affine bundles, the mentioned construction
of the structure group calculus gives us an intrinsic bicovariant *-calculus $\aC$ over 
$\aG$. It is determined by the system of affine maps 
$\varrho_D, \chi_E\colon\aA\rightarrow \ahorP$.

Before passing to the explicit construction of the affine calculus $\aC$,
let us write down some algebraic relations involving the maps
$\varrho_D$ and $\chi_E$. We have already introduced the graded
commutant $\zhP$. It is obviously included in the `affine' graded commutant
$$\azhP\leftrightarrow\zhP\otimes\SV$$ 
which is a right $\aA$-module, in a natural manner (the module
structure is constructed with the help of the affine translation map). Explicitly, the right
$\aA$-module structure on $\azhP$ is specified by
\begin{gather*}
(\zeta\otimes\vartheta)\circ a=(\zeta{\circ} a^{(1)})\otimes(\vartheta{\circ} a^{(2)})\\
\xi\circ\vartheta=\sum_{i\geq 0}\sum_{\alpha,\beta\in S[i]}(-1)^i \mathrm{I}\{t_{\beta i}\}
(\xi{\circ}c_{\beta\alpha}^i) \vartheta_{\alpha i},
\end{gather*}
where $\xi\in \azhP$, $a\in\cal{A}$ and $\vartheta\in\SV$. 

The maps $\varrho_D,\chi_E\colon\aA\rightarrow\ahorP$ take their
values from $\azhP$. Furthermore, according to the general theory \cite{d-diff}, the 
maps $\varrho_D$ and $\chi_E$, restricted to $\ker\bigl(\e\colon\aA
\rightarrow\Bbb{C}\bigr)$, intertwine the
right $\aA$-multiplication and the right $\aA$-module structure on $\azhP$.
Explicitly, we have
\begin{gather}
\varrho_D(\xi a)=\varrho_D(\xi)\circ a=\lt{a}\varrho_D(\xi)\rt{a}\qquad\chi_E(\xi a)
=\chi_E(\xi)\circ a\label{aDE1}\\
\begin{aligned}\label{aDE3}
\varrho_D(\xi\vartheta)&=\varrho_D(\xi)\circ\vartheta=
\sum_{i\geq 0}\sum_{\alpha,\beta\in S[i]}(-1)^i \mathrm{I}\{t_{\beta i}\}\bigl\{
\varrho_D(\xi){\circ}c_{\beta\alpha}^i\bigr\} \vartheta_{\alpha i} \\
\chi_E(\xi \vartheta)&=\chi_E(\xi)\circ\vartheta=
\sum_{i\geq 0}\sum_{\alpha, \beta\in S[i]}(-1)^i \mathrm{I}\{t_{\beta i}\}
\bigl\{\chi_E(\xi){\circ}c_{\beta\alpha}^i \bigr\}\vartheta_{\alpha i},
\end{aligned}
\end{gather}
where the notation is the same as in formula \eqref{taP1}, with $a\in\cal{A}$ and
$\vartheta\in\SV$ while
$\xi\in\aA$ satisfies $\e(\xi)=0$. In particular, for $\vartheta=\theta\in\V$ 
formulae \eqref{aDE3} reduce to 
\begin{equation}\label{aDE2}
\begin{aligned}
\varrho_D(\xi \theta)&=\varrho_D(\xi)\circ\theta=\varrho_D(\xi)\theta-
\sum_k\theta_k\{\varrho_D(\xi)\circ c_k\}\\
\chi_E(\xi \theta)&=\chi_E(\xi)\circ\theta=\chi_E(\xi)\theta-
\sum_k\theta_k\{\chi_E(\xi)\circ c_k\}.
\end{aligned}
\end{equation} 

The above expressions allow us to construct explicitly the calculus $\aC$, 
in terms of the calculus $\Gamma$ and the algebraic relations in 
$\ahorP$. In particular, we see that the space $\aCi$ contains in a 
natural manner $\Gamma_{\inv}$ and the image of $\V$. Moreover 
$\aCi$ is built from $\Gamma_{\inv}$ and $\V$, by the right action of $\SV$. 

\smallskip
******\par
The presented construction of the calculus $\aC$ can be easily abstracted from 
the context of quantum principal bundles. We are going to make a small digression from
our main theme, and to present an abstract version of the construction of $\aC$. 

Let us consider a *-algebra $\cal{H}$, 
equipped with the right $\cal{A}$-module structure $\circ$ and a right 
counital $\cal{A}$-comodule structure 
$\Delta\colon\cal{H}\rightarrow\cal{H}\otimes\cal{A}$, such that 
the standard compatibility properties
\begin{gather*}
(\varphi\circ a)^*=\varphi\circ \k(a)^* \qquad\Delta*=(*\otimes *)\Delta\\
(\varphi\psi)\circ a=(\varphi\circ a^{(1)})(\psi\circ a^{(2)})\\
\Delta(\varphi\circ a)=\sum_k\bigl(\varphi_k\circ a^{(2)}\bigr)\otimes\k(a^{(1)})c_k a^{(3)}
\end{gather*}
hold, where $\Sum_k\varphi_k\otimes c_k=\Delta(\varphi)$. 

Let us now introduce, starting from the algebra $\cal{H}$, a *-algebra 
$\widetilde{\cal{H}}=\cal{H}\otimes\SV$, equipped with the standard 
cross-product structure. The maps $\circ$ and $\Delta$ are naturally
extendible to the right $\aA$-module structure $\circ$ on $\widetilde{\cal{H}}$
and the right comodule structure $\widehat{\Delta}\colon
\widetilde{\cal{H}}\rightarrow\widetilde{\cal{H}}\otimes\aA$, such that 
the above compatibility properties are preserved. Explicitly, the map 
$\widehat{\Delta}$ is constructed by combining $\Delta$ and $\widehat{\v}$. 
The extension of $\circ$ is constructed in two steps. At first, right $\cal{A}$-module
$\circ$-structures on $\cal{H}$ and $\SV$ naturally combine to a right 
$\cal{A}$-module structure on $\wH$. Secondly, the formula
\begin{equation}
\xi\circ\vartheta=
\sum_{i\geq 0}\sum_{\alpha,\beta\in S[i]}(-1)^i \mathrm{I}\{t_{\beta i}\}
(\xi{\circ}c_{\beta\alpha}^i) \vartheta_{\alpha i},
\end{equation}
uniquely and consistently defines an extension of $\circ$ to a right 
$\aA$-module structure on $\wH$. Here $\vartheta\in\SV$ and the notation is 
the same as in formula \eqref{taP1}. 

Let us assume that a map $\pi\colon\cal{A}\rightarrow\cal{H}$ is given, 
such that 
\begin{gather*}
\Delta\pi=(\pi\otimes\id)\ad\qquad\pi(a)^*=-\pi[\k(a)^*]\\
\pi(ab)=\pi(a){\circ}b+\e(a)\pi(b),\quad \Rightarrow \quad \pi(1)=0. 
\end{gather*}
This map naturally defines a bicovariant *-calculus $\Gamma$ over $G$. It is based
on the right $\cal{A}$-ideal $\cal{R}=\ker(\e)\cap\ker(\pi)$. 

Furthermore we shall assume that a map $\zeta\colon\V\rightarrow\cal{H}$ is given, 
satisfying 
\begin{equation*}
\Delta\zeta=(\zeta\otimes\id)\v, \qquad\zeta *=*\zeta,\qquad 
\zeta(\theta{\circ}a)=\zeta(\theta){\circ}a.  
\end{equation*}
We are going to construct, with the help of $\zeta$, 
a natural extension $\wpi$ of $\pi$, which is defined on $\aA$ and takes the
values from $\wH$. Firstly, let us define 
$(\wpi\restr\V)\colon\V\rightarrow\wH$ by
\begin{equation}\label{wpi-restr}
\wpi(\theta)=\zeta(\theta)-\sum_k\theta_k\pi(c_k),
\end{equation}
where $\Sum_k\theta_k\otimes c_k=\v(\theta)$. The full map 
$\wpi\colon\aA\rightarrow\wH$ is
defined by the following construction:

\begin{pro} The formula 
\begin{equation}\label{defwpi}
\wpi(q\vartheta)=
\sum_{i\geq 0}\sum_{\alpha,\beta\in S[i]}(-1)^i \mathrm{I}\{t_{\beta i}\}\bigl\{
\wpi(q){\circ}c_{\beta\alpha}^i \bigr\}\vartheta_{\alpha i}
\end{equation}
consistently defines a map $\wpi\colon\aA\rightarrow\wH$, extending the map
$\pi\colon\cal{A}\rightarrow\cal{H}$ and previously defined restriction $\wpi\restr\V$. 
Here $q\in\ker\bigl(\e:\aA\rightarrow \Bbb{C}\bigr)$ and $\vartheta\in\SV$. The following identities hold:
\begin{gather}
\widehat{\Delta}\wpi=(\wpi\otimes\id)\adaG\qquad\wpi(\xi)^*=-\wpi[\k(\xi)^*]\label{wpi-*-cov}\\
\wpi(\xi\psi)=\wpi(\xi){\circ}\psi+\e(\xi)\wpi(\psi). \label{wpi-circ}
\end{gather}
In particular, $\wpi$ defines naturally a bicovariant *-calculus $\aC$ over $\aG$. The space
$\aCi=\aA/\ker(\wpi)\leftrightarrow\im(\wpi)$ is a right $\aA$-submodule of $\wH$. As a 
right $\SV$-module, it is generated by $\Gamma_{\inv}$ and $\wpi(\V)$. 
\end{pro}

\begin{pf} Let us first check the consistency of \eqref{defwpi}. One method is to proceed
directly, by checking the compatibility of \eqref{defwpi} with the commutation 
relations (twisted commutation between $\cal{A}$ and $\SV$, and the kernel of 
the braided symmetrizator map defining $\SV$) defining the algebra $\aA$. Here we shall 
sketch a different method, by constructing first the associated bicovariant algebras. 

Let $\mho\leftrightarrow\cal{A}\otimes \cal{H}$ be a bicovariant *-algebra associated to 
$(\cal{H},\Delta, \circ)$. The product in $\mho$ is given by the standard cross-product rule. 
The maps $\ell_\mho=(\phi\otimes\id)\colon\mho\rightarrow \cal{A}\otimes\mho$ and 
$\wp_\mho=(\id\otimes\Delta)\colon\mho\rightarrow\mho\otimes\cal{A}$ are the corresponding
left/right actions of $G$. Let us observe that the formula 
$$ da=a^{(1)}\otimes \pi(a^{(2)})$$
defines a hermitian differential $d\colon\cal{A}\rightarrow \mho$ which is left/right-covariant, in a 
natural way. Let $\Gamma$ be the $\cal{A}$-bimodule defined as 
$$\Gamma=\Bigl\{\Sum ad(b)\Bigm\vert a,b\in\cal{A}\Bigr\}.$$
By construction we see that $\Gamma$, equipped with $d$, is a bicovariant *-calculus over $G$. 

In a similar way, we can construct a bicovariant *-algebra 
$\amalg\leftrightarrow\aA\otimes\widetilde{\cal{H}}$, equipped with the left/right actions
$\ell_\amalg\colon\amalg\rightarrow \aA\otimes \amalg$ and $\wp_\amalg\colon\amalg\rightarrow
\amalg\otimes\aA$. The algebra $\mho$ is a $\aG$-covariant *-subalgebra of $\amalg$. 
Furthermore, there exists a unique differential $d\colon\aA\rightarrow\amalg$ extending 
the map $d\colon\cal{A}\rightarrow\mho$ and satisfying 
$$ d(\theta)=\zeta(\theta)\qquad\forall\theta\in\V. $$
It is easy to see that the extended $d$ is also hermitian,  and left/right $\aG$-covariant. Let
$\aC$ be the corresponding bicovariant *-calculus over $\aG$. We have $\aCi\subseteq\widetilde{\cal{H}}$ 
and $\Gamma_{\inv}\subseteq\cal{H}$. Now we can define the extension 
$\wpi\colon\aA\rightarrow\aCi$ of $\pi\colon\cal{A}\rightarrow\Gamma_{\inv}$ by the formula
$$ \wpi(\psi)=\k(\psi^{(1)})d(\psi^{(2)}). $$
It is easy to see that \eqref{wpi-restr} and \eqref{defwpi} holds. Finally, properties 
\eqref{wpi-*-cov} and \eqref{wpi-circ} are the standard properties of the quantum germ maps. 
\end{pf}

Let us observe that if the group $G$ is {\it $\Gamma$-connected} (in the sense that 
only the scalar elements of $\cal{A}$ are annihilated by the differential $d$) and if there are no
$\v$-invariants in $\V$, then the 
restriction map $\wpi\colon\V\rightarrow \aCi$ will always be injective, independently of $\zeta$. 
In what follows we shall assume that the injectivity property holds. The grading on $\SV$ enables 
us to introduce a natural filtering on $\aCi$, compatible with the right
$\SV$-module structure. For $k\geq0 $ let us define the spaces 
$$ \aCi^k=\wpi\bigl[\SV^{k+1}\bigr]+\Gamma_{\inv}{\circ}\SV^k. $$
In general, these spaces do not form a direct sum, however they span the whole $\aCi$ and
we can introduce a filtering $\bigl\{{}_m\!\aCi\bigr\}$ of $\aCi$ by defining 
${}_m\!\aCi=\Sum_{k\leq m}\aCi^k$. We have $\aCi^0=\Gamma_{\inv}\oplus\V$. 

\medskip
Going back to our main considerations, we see that the calculus $\aC$ is the minimal 
calculus compatible with the system of maps $\varrho_D$ and $\chi_E$. 
Let us also observe that the appearance of the non-trivial filtering in $\aCi$ is a 
purely quantum phenomena. By construction, the curvature/transition maps are 
factorizable, and can be interpreted as $\chi_E,\varrho_D\colon\aCi \rightarrow\ahorP$.
\begin{lem} We have
\begin{equation}
\chi_E,\varrho_D(\aCi^k)\subseteq
\horP^2\otimes\SV^k+\horP^2\otimes\SV^{k+1}\qquad\forall k\in\Bbb{N}\cup\{0\}. 
\end{equation}
In particular, if $|k-l|\geq 2$ then $\aCi^k\cup\aCi^l=\{0\}$. 
\end{lem}
\begin{pf} The above inclusions follow by applying formulae \eqref{ext-r} and \eqref{ext-k} 
together with property \eqref{aDE3}. 
\end{pf}

The higher-order spaces $\aCi^k$ will be trivial only in some very special cases (including 
the classical case). 
 
\subsection{The Global Calculus}

We can now continue with the construction of the global calculus on the 
quantum affine bundle $\aP$, following \cite{d-diff}. The construction builds, 
in an intrinsic way, a graded-differential *-algebra $\aWP$, starting from
 *-algebra $\ahorP$, together with the system of maps
$D,E\colon\ahorP\rightarrow\ahorP$ and the calculus $\aC$.
Every covariant derivative $D$ induces an identification
\begin{equation}
\aWP\leftrightarrow\ahorP\otimes\aCi^\wedge=\vh(\aP)
\end{equation}
of graded vector spaces. The graded *-algebra structure of $\vh(\aP)$ is
given by the twisted product between $\ahorP$ and $\aCi^\wedge$. In
accordance with \cite{d-diff}, there exists
a natural bijection between  affine
covariant derivative maps $D=D_\omega$ and regular connections
$\omega=\omega_D\colon\aCi\rightarrow\aWP$ on $\aP$.

The elements of $\aWP$ play the role of
differential forms on the quantum affine bundle $\aP$. There exists a
natural homomorphism
$\wFh\colon\aWP\rightarrow\aWP\grten\aCi^\wedge$ of
graded-differential *-algebras, extending the map $\Fah$. We have
$$ \ahorP=\wFh^{-1}\Bigl\{\aWP\otimes\aA\Bigr\}. $$

It is important to mention that our constructions of the minimal calculi $\Gamma$ 
and $\aC$, as well as the global differential calculi on $P$ and $\aP$, are based on
the complete space of bundle derivatives. These constructions are applicable, without
any change, to an arbitrary affine subspace $\cal{L}\subseteq \DP$. In general, such 
a modified construction will give us a simpler calculus, which will be not compatible
with all possible maps from $\DP$. However, in various concrete situations even {\it a 
single derivative} $\{D\}=\cal{L}$ will generate the complete calculus. In classical 
geometry, for example, such special connections are those having {\it the maximal} 
infinitezimal holonomy (the corresponding calculus will be the classical calculus on 
the structure group and the bundle). 

\section{Concluding Examples and Observations}
\subsection{Frame Structures and Affine Bundles}

In this subsection we shall briefly analyze relations between
constructed affine bundles and
general frame structures \cite{d-frm2} on quantum principal bundles.
Let $\V^\wedge$ be the $\tau$-exterior algebra associated to $\V$. This is
a quadratic algebra given by the relations
$$ \im(I+\tau)\subseteq\V\otimes\V, $$
and therefore it is natural to assume that $\ker(I+\tau)\neq\{0\}$. The algebra
$\V^\wedge$ is equipped with the induced $\circ$-structure.

Let us consider a quantum principal $G$-bundle $P=(\cal{B},i,F)$. Let
$\hor_P$ be a graded *-algebra given by
$$ \hor_P\leftrightarrow\cal{B}\otimes\V^\wedge$$
at the level of graded vector spaces,
with the *-algebra structure specified by equalities
\begin{align*}
(q\otimes\vartheta)(b\otimes\eta)&=\sum_kqb_k\otimes(\vartheta{\circ} c_k)
\eta\\
(b\otimes\vartheta)^*&=\sum_k b_k^*\otimes(\vartheta^*{\circ} c_k^*),
\end{align*}
where $\Sum_k b_k\otimes c_k=F(b)$.

The elements of $\hor_P$ are interpretable as `frame-type' horizontal
forms on the bundle $P$. There exists
a unique homomorphism
$\Fh\colon\hor_P\rightarrow\hor_P\otimes\cal{A}$ extending both
$F$ and $\v$. This map is coassociative and hermitian.

The constructed algebra $\hor_P$ is in the roots of the
abstract definition of quantum frame structures \cite{d-frm2}. Quite precisely,
a quantum frame structure is given by a $\Fh$-covariant hermitian
first-order antiderivation $\nabla$ on $\hor_P$ such that its restriction
$\dM\colon\Omega_M\rightarrow\Omega_M$ on
the $\Fh$-invariants is a differential having the property
that $\cal{V}=\Omega_M^0$ and $\dM(\cal{V})$ generate the module
$\Omega^1_M$ (and hence the whole $\Omega_M$), and such that $\nabla(\V)=\{0\}$.
The map $\nabla$ plays the role of the Levi-Civita connection, and the elements of
$\V^\wedge$ are counterparts of the `coordinate forms' on the frame bundle.

In the case of the frame structures, there exists an intrinsic choice of the 
translational part $\lambda\in\Tcon^1$. It is simply given by 
\begin{equation}
\lambda(\theta)=\theta_\wedge, 
\end{equation}
where $\theta\in\V$ and we have denoted by $\theta_\wedge$ the elements of $\V$, 
understood as coordinate 1-forms (the elements of $\V^\wedge$). 

For arbitrary affine connections constructed with the help of this `reinterpretation'
map, Proposition \ref{pro-affcurv} allows us to view the torsion tensor as a part of the
affine curvature. Indeed, the torsion of a bundle derivative $D$ is defined as 
$T(\theta)=D(\theta_\wedge)$, as in classical geometry. Hence we have 
\begin{equation}
\varrho_D(\theta)=-\sum_k\vartheta_k\varrho_D(c_k)-T(\theta),
\end{equation} 
which is in a complete agreement with classical geometry \cite{KN}. 

As a concrete illustration, let us consider the quantum Hopf fibration $P$. This is
a quantum $\U(1)$-bundle over a quantum $2$-sphere \cite{p}, defining a spin strucrure. 
The total space of the bundle $P$ is the quantum $\SU(2)$-group \cite{w-su2}.  

The Hopf *-algebra $\cal{A}$ describing $G$ is generated by a single unitary element $U$,
corresponding to the canonical inclusion of the unit circle in the complex plane $\Bbb{C}$.
It is equipped with the group law $$\phi(U)=U\otimes U.$$

By definition \cite{w-su2}, the *-algebra $\cal{B}$  is
generated by elements $\alpha$ and $\gamma$ together with the relations 
\begin{gather*}
\alpha\alpha^*+\mu^2\gamma\gamma^*=1\qquad\alpha^*\alpha+\gamma^*\gamma=1\\
\alpha\gamma=\mu\gamma\alpha\qquad
\alpha\gamma^*=\mu\gamma^*\alpha\qquad
\gamma\gamma^*=\gamma^*\gamma,
\end{gather*}
where $\mu\in[-1,1]\setminus\{0\}$. 
The fundamental representation is defined by 
\begin{equation*}
u=\begin{pmatrix}
\alpha&-\mu\gamma^*\\
\gamma&\phantom{-\mu}\alpha^*
\end{pmatrix}
\end{equation*} 
and the above relations are equivalent to the unitarity of this matrix. The coproduct is 
given by the standard matrix rule
$$ \phi(u_{ij})=\sum_k u_{ik}\otimes u_{kj}. $$

Let $\Phi$ be the $3$-dimensional calculus over $P$ constructed in \cite{w-su2}. This
calculus is left-covariant and *-covariant. 
The elements 
$$\eta_3=\pi(\alpha-\alpha^*)\qquad
\eta_+=\pi(\gamma)\qquad
\eta_-=\pi(\gamma^*)$$
span the space $\Phi_{\inv}$. The associated right ${\cal B}$-module structure 
$\circ:\Phi_{\inv}\otimes\cal{B}\rightarrow \Phi_{\inv}$ is specified by 
\begin{equation*} 
\begin{aligned}
\mu^2\eta_3\circ\alpha&=\eta_3\\
\mu\eta_\pm\circ\alpha&=\eta_\pm
\end{aligned}\qquad
\begin{aligned}
\eta_3\circ\alpha^*&=\mu^2\eta_3\\
\eta_\pm\circ\alpha^*&=\mu\eta_\pm
\end{aligned}
\end{equation*}
with $\Phi_{\inv}{\circ}\gamma=\Psi_{\inv}{\circ}\gamma^*=\{0\}$.  Here we have used the
same symbols for the group entities operating on $G$ and $P$. 

We shall identify $G=\U(1)$ with a classical subgroup of $P$ consisting of diagonal 
matrices (if $\mu\neq -1,1$ then $G$ is precisely the classical part of $P$). In other words, 
$$\cal{A}\leftrightarrow\cal{B}/\mathrm{gen}(\gamma, \gamma^*).$$
More precisely $U\leftrightarrow [\alpha]$ and $U^*\leftrightarrow [\alpha^*]$. 
We see that the above $\circ$-structure on $\Phi_{\inv}$ is naturally projectable down to a
right $\cal{A}$-module structure on the same space (and will be denoted by the same symbol). 

Let us define $\V$ as the space spanned by $\eta_\pm$. In what follows, $\V$ will be
equipped with the induced $\circ$ and $*$-structures. 
Furthermore, the map $\chi\colon\V\rightarrow\V\otimes \cal{A}$ is defined as {\it the adjoint} action 
of $G$. Explicitly,
\begin{equation}
\chi(\eta_+)=\eta_+\otimes U^2\qquad\chi(\eta_-)=\eta_-\otimes U^{-2}. 
\end{equation}
 
The canonical braid operator $\tau\colon\V^{\otimes 2}\rightarrow\V^{\otimes 2}$ 
has a simple form
\begin{equation}
\tau=\begin{pmatrix}1/\mu^2& 0 & 0 & 0\\
                         0 & 0 & 1/\mu^2 & 0\\
                         0 & \mu^2 & 0 & 0\\
                         0 & 0 & 0 & \mu^2 \end{pmatrix}
\end{equation}
in the basis $\eta_\pm$. It follows that the $\tau$-exterior algebra is given by
the relations 
\begin{equation}
\eta_\pm^2=0\qquad\quad\eta_+\eta_-=-\mu^2\eta_-\eta_+. 
\end{equation}

These relations are just a subset of the full set of relations defining the canonical 
higher-order calculus over $P$. This calculus is given by the universal differential
envelope $\Phi^\wedge$ of $\Phi$ (\cite{w-su2},\cite{d1}--Appendix B). The remaining 
relations are 
\begin{equation}
\eta_3^2=0, \qquad\quad \eta_3\eta_\pm=\mu^{\mp 4}\eta_\pm\eta_3.
\end{equation}

The algebra $\horP$ is therefore the subalgebra of $\Phi^\wedge$ generated by 
$\cal{B}$ and $\eta_+,\eta_-$. By construction, there exists a natural projection homomorphism 
$p_{h\!or}\colon \Phi^\wedge\rightarrow \horP$ defined by annihilating the vertical
part $\eta_3$. Composing this projection with the differential 
$d\colon\Phi^\wedge\rightarrow\Phi^\wedge$ we obtain a first-order 
*-antiderivation $\nabla\colon\horP\rightarrow\horP$. This map gives us a canonical
frame structure, and it coincides with the covariant derivative of the canonical 
regular connection introduced in \cite{d2}. It corresponds to the standard 
Levi-Civita connection on the $2$-sphere (lifted to the corresponding spin bundle). 

As already mentioned, in accordance with the general theory \cite{d-frm2}, 
the map $\nabla$ induces a canonical differential calculus $\Gamma$ on the structure group 
$G=\U(1)$. In this case the calculus is $1$-dimensional, spanned by
the element $\zeta=\pi(U-U^*)$. This calculus is non-classical, because 
the $\circ$-structure on $\Gamma_{\inv}$ is not trivial. Actually $\Gamma$ is the 
projection of $\Phi$ on $G$ and we have 
$$\zeta\circ U=\frac{1}{\mu^2}\zeta. $$
Furthermore, the curvature map is given by 
\begin{gather}
\rho_\nabla(U^n)=\frac{1-\mu^{-2n}}{1-\mu^{-2\phantom{n}}}\rho_\nabla(U)\qquad\forall 
n\in\Bbb{Z}\\
\rho_\nabla(U)=\mu\eta_-\eta_+\qquad \rho_\nabla(U^{-1})=\mu\eta_+\eta_-.
\end{gather}

Now let us construct the affine extension $\aG$ of $G$. We shall use the symbols
$\xi\leftrightarrow\eta_+$ and $\xi^*\leftrightarrow\mu\eta_-$ when referring to 
the braided-symmetric algebra $\SV$ build over $\V$. It is easy to see that 
the algebra $\SV$ is determined by a single relation
\begin{equation}
\xi\xi^*=\mu^2\xi^*\xi. 
\end{equation}

The corresponding Hopf *-algebra $\aA$ is determined by the unitary $U$
generating $\cal{A}$, generators $\{\xi,\xi^*\}$ defining $\SV$, and additional relations
\begin{equation}
\xi U=\frac{1}{\mu}U\xi\quad \Leftrightarrow\quad \xi^* U=\frac{1}{\mu}U\xi^*,
\end{equation} 
together with the coproduct specifications
\begin{equation}
\phi(U)=U\otimes U\qquad \phi(\xi)=1\otimes\xi+\xi\otimes U^2. 
\end{equation}

Geometrically, the constructed quantum group is the two-fold covering of the quantum 
$\E(2)$ group. The bundle $P$ is interpretable as the spin $2$-fold covering of 
the orthonormal frame bundle of the quantum sphere. 

Let us naturally extend the Levi-Civita connection to the affine level, as 
explained at the beginning of this subsection. We are going to construct the calculus 
$\aC$, associated to the extended Levi-Civita connection $\nabla$. Let us assume that 
$\mu\in(-1,1)\setminus\{0\}$. 

\begin{pro}We have 
\begin{equation}\label{r-xi}
\varrho_\nabla(\xi^p\xi^{*q})=c_{pq}w \otimes \xi^p\xi^{*q} \qquad\forall p,q\quad |p|+|q|> 0, 
\end{equation}
where 
$$
w=\frac{\mu^3}{1-\mu^2}\eta_-\eta_+
$$
and
$$
c_{pq}=\prod_{i=1}^p(1-\frac{1}{\mu^{2i+2}})\prod_{j=1}^q(1-\mu^{2j+2}).
$$
 
In particular, it follows that the space $\aCi$ is infinite-dimensional, 
with a natural basis given by the elements $\{e_{pq}\mid p,q\geq 0\}$, where 
$$
e_{pq}=\begin{cases}\zeta& \text{for $p=q=0$},\\
\wpi(\xi^p\xi^{*q})&\text{otherwise.}\end{cases}\qquad 
$$
The vector $\zeta$ is cyclic and separating for the $\circ$-action of $\SV$ on $\aCi$. 
\end{pro}

\begin{pf}
At first, let us evaluate the curvature on the generators $\xi,\xi^*$. 
In accordance with \eqref{ext-r} and the torsionless property of $\nabla$ we have
\begin{gather*}
\rho_\nabla(\xi)=-\xi\rho_\nabla(U^2)=(1-\frac{1}{\mu^4})w \otimes \xi\\
\rho_\nabla(\xi^*)=(1-\mu^4)w\otimes \xi^*. 
\end{gather*}
This proves the cases $(p,q)=(0,1)/(1,0)$. The general expression \eqref{r-xi}
follows by induction on $p,q$ and applying the first of \eqref{aDE2}. To complete
the proof, it is sufficient to observe that the image $\rho_\nabla(\aA)$ is spanned
by $\rho_\nabla(U-U^*)$ and the elements \eqref{r-xi}. 
\end{pf}

If $\mu=1$ we are in the framework of classical geometry, and all objects are classical. 
In particular $P$ is the classical Hopf fibration, and $\aC$ is the standard $3$-dimensional 
calculus over the spin covering $\aG$ of $\E(2)$. If $\mu=-1$ the bundle $P$ will
be quantum, however the base is the classical $2$-sphere. Furthermore, the calculus $\aC$ 
as well as the group $\aG$ are both non-classical. The space $\aCi$ is $3$-dimensional
(spanned by $\zeta$, $\wpi(\xi)$ and  $\wpi(\xi^*)$) with $\aCi^k=\{0\}$ for $k\neq 0$,
and the $\circ$-structure is given by 
$$\{\xi,\xi^*\}\circ U=-\{\xi,\xi^*\}\qquad \{\xi,\xi^*\}\circ \V=\{0\}. $$
The braid operator $\tau$ is the standard transposition. 

The presented calculations are easily incorporable in the context of arbitrary quantum frame 
$\U(1)$-bundles. Geometrically, such structures are interpretable as quantum 
Riemann surfaces. However, for general surfaces the Levi-Civita connection $\nabla$ will
generate a multi-dimensional calculus over the structure group $G=\U(1)$. This 
phenomena will occur if and only if $\rho_\nabla$ is not an eigenvector for
the right $\cal{A}$-module structure $\circ$ on $\zhP$. More precisely, 

\begin{lem} Assume that the curvature of $\nabla$ is non-zero. The space $\Gamma_{\inv}$ is 
naturally identificable with $\im[\rho_\nabla]$. The image of $\rho_\nabla$ is 
the $\circ$-submodule of $\zhP$ generated by $\rho_\nabla(U)$. \qed
\end{lem}

\subsection{Affine Structures as Frame Structures}

As we have seen, every bundle derivative acting on a frame bundle
can be naturally extended to the affine level, so that the curvature of the
extended connection appears as a combination of the original curvature 
(involving translational coordinates) and the torsion operator. Let us now consider 
the opposite problem---to construct the frame structure from a given affine 
connection. 

Our starting point is an arbitrary abstract horizontal algebra $\horP$ 
as defined in Section~\ref{sec:ahor}, and we shall also assume that 
$\WM$ is generated, as a differential algebra, by $\WM^0=\cal{V}$. 
Let us consider a map $\lambda\in\Tcon^1$. The following relations hold
\begin{gather} 
-\lambda(\eta)\lambda(\theta)=\sum_k\lambda(\theta_k)\lambda(\eta{\circ}c_k)\\
\lambda(\theta)b=\sum_\alpha b_\alpha\lambda(\theta{\circ} d_\alpha), 
\end{gather} 
as it follows from \eqref{v-L}--\eqref{hor-L}. Here $\v(\theta)=\Sum_k\theta_k\otimes c_k$ and
$\Sum_\alpha b_\alpha\otimes d_\alpha=F(b)$. This implies that there exists a unique homomorphism 
$h_\lambda\colon\cal{B}\otimes\V^\wedge\rightarrow\horP$,
reducing to the identity on $\cal{B}$, and extending $\lambda$. The space
$\cal{B}\otimes\V^\wedge$ is equipped with the standard cross-product structure. Moreover,
the map $h_\lambda$ intertwines the corresponding actions of $G$.

\begin{defn} We shall say that $\lambda\in\Tcon^1$ is {\it regular} iff 
the map $h_\lambda$ is bijective.
\end{defn}

In other words, we can identify $\cal{B}\otimes\V^\wedge$ and $\horP$
with the help of $h_\lambda$. 
If in addition there is a bundle derivative $\nabla\in\DP$ such that 
$\im(\lambda)\subseteq\ker(\nabla)$ then $\nabla$, together with 
the identification $\cal{B}\otimes\V^\wedge\leftrightarrow\horP$ 
determines a frame structure on the bundle $P$. On the other hand, 
$\nabla$ and $\lambda$ determine an extended bundle derivative
$\nabla_*\colon \ahorP\rightarrow\ahorP$. The map $\nabla_*$ 
contains the full information about the constructed frame structure
and by construction it is the Levi-Civita connection associated to
this frame structure. 

In the context of the frame structures, the space $\ahorP$ is decomposed as follows:
$$
\ahorP\leftrightarrow \cal{B}\otimes\negthickspace\negthickspace\negthickspace\underbrace{\V^\wedge
\otimes\SV}_{\text{\normalsize\it q-plane $\Omega(\V,\tau)$}}\negthickspace\negthickspace\negthickspace.
$$
It is interesting to observe that the marked graded *-subalgebra is actually the 
canonical differential calculus on the quantum plane built over $(\V,\tau)$. The
Levi-Civita connection $\nabla_*$ reduces to the standard differential in $\Omega(\V,\tau)$. 

\subsection{Analogy With Horizontal-Vertical Decompositions}

There exists an interesting formal analogy between the formalism of
affine bundles $\aP$ and differential calculus on general quantum
principal bundles $P$. In this analogy, $\V$ corresponds to the space
$\Gamma_{\inv}$. The analogy is more complete if we
assume that the higher-order calculus on $G$ is based on the 
braided exterior algebra $\Gamma_{\inv}^\vee$, because 
$\tau\leftrightarrow -\sigma$ and $\SV\leftrightarrow\Gamma_{\inv}^\vee$. 
Then the algebra $\ahorP$ corresponds to the algebra $\vh(P)$ of
horizontally-vertically decomposed differential forms.

Let us now assume that $\gaP=(\gaB,i,\gFh)$ is an arbitrary quantum principal
$\aG$-bundle over a quantum space $M\leftrightarrow\cal{V}$. We can introduce a filtration on
$\gaB$, by the formula
$$
\gaB_k=\gFh^{-1}\Bigl\{\gaB\otimes\aA_k\Bigr\},\qquad\aA_k=\Bigl\{
\sum a\vartheta\Bigm\vert a\in\cal{A}, \vartheta\in \SV^j, j\leq k\Bigr\}.
$$
In particular, let us consider the *-subalgebra $\cal{B}=\gaB_0$. We see
that this algebra is the analog of horizontal forms in
the theory of differential calculus. In other words, we have
$$
\gFh(\cal{B})\subseteq\cal{B}\otimes\cal{A}.
$$
Let us observe that $i(\cal{V})\subseteq\cal{B}$ and that $P=(\cal{B},i,F)$ is
a quantum principal $G$-bundle over $M$. Here $F$ is the restricted coaction map. 

\begin{defn} A {\it translaton} on $\aP$ is every hermitian linear map
$\xi\colon\V\rightarrow\gaB$ satisfying
\begin{equation}
\gFh\xi(\theta)=\sum_k\theta_k\otimes c_k +1\otimes\theta,
\end{equation}
where $\v(\theta)=\Sum_k\theta_k\otimes c_k$.
\end{defn}

The translatons are the affine-geometrical analogs of connections, as defined in
the general theory. Every affine bundle $\aP$ admits at least one
translaton. The proof is similar as in the connection existence
theorem \cite{d2}. It is based on the representation theory \cite{w2}
of compact matrix quantum groups. The translatons of $\aP$
form a real affine space. Furthermore, it is natural to
formulate

\begin{defn} A translaton $\xi$ is {\it regular} iff
\begin{equation}
\xi(\theta)\varphi=\sum_j\varphi_j\xi(\theta\circ d_j),
\end{equation}
where $\varphi\in\cal{B}$, and $\Sum_j\varphi_j\otimes d_j=\gFh(\varphi)$.
We shall say that $\xi$ is {\it multiplicative}, iff it extends
(necessarily uniquely)
to a unital homomorphism $\xi\colon\SV\rightarrow\gaB$. Equivalently,
$\xi^\otimes[\ker(Y)]=\{0\}$.
\end{defn}

The lack of multiplicativity of $\xi$ is measured by the
values of $\xi^\otimes$ on the generators of the ideal $\ker(Y)$. Similarly,
the lack of regularity of $\xi$ is measured by the values
of the operator
$\ell^\xi\colon\V\otimes\cal{B}\rightarrow\gaB$, defined by
$$
\ell^\xi(\theta,\varphi)=\xi(\theta)\varphi-
\sum_j\varphi_j\xi(\theta\circ d_j).
$$
It is worth noticing that
\begin{align}
\gFh\ell^\xi(\theta,\varphi)&=
\sum_{kj}\ell^\xi(\theta_k,\varphi_j)\otimes c_kd_j\\
[\ell^\xi(\theta,\varphi)]^*&=-\sum_j\xi(\theta^*\circ\k(d_j)^*,\varphi_j^*).
\end{align}
The operator $\ell^\xi$ always
takes its values from the algebra $\cal{B}$.
Let us fix a grade-preserving *-invariant and $\v$-covariant splitting
$$
\V^\otimes=\ker(Y)\oplus\SV.
$$

This splitting induces a map
$m_\xi\colon\cal{B}\otimes\SV\rightarrow\gaB$.
It turns out that this map is bijective, and intertwines the
corresponding actions of $G$. Moreover, the map $m_\xi$ will
be a *-homomorphism iff $\xi$ is regular and multiplicative.

In any case, $m_\xi$ induces an intrinsic graded *-isomorphism
$$
\cal{B}\otimes\SV\leftrightarrow\gr\bigl\{\gaP\bigr\},
$$
which is independent of the choice of the translaton $\xi$.

Therefore, the affine bundles introduced in a constructive way
in Section~2 can be characterized as principal $\aG$-bundles admitting
regular and multiplicative translatons.

\end{document}